\documentclass[3p,10pt]{elsarticle}

\usepackage{latexsym,epsf,amssymb,amsthm,amsfonts,multicol,makeidx}
\usepackage{color}
\usepackage{colortbl,xcolor}
\usepackage{cancel}

\definecolor{marron}{rgb}{0.5,0.3,0}

\usepackage{epsfig,amssymb}
\usepackage{amsmath}
\usepackage{dsfont,upgreek}
\usepackage{relsize}

\usepackage{multicol}
\usepackage{multirow}

\usepackage{tikz}
\usepackage{pgfplots}
\usepackage{mathrsfs}
\usetikzlibrary{arrows}
\usetikzlibrary[patterns]

\def\doubleone{ 1 \kern-.224em \hbox{\rm l}}

\def\doubleone{ 1 \kern-.224em \hbox{\rm l}}

\newcommand{\bq}{\begin{equation}}
\newcommand{\eq}{\end{equation}}

\newcommand{\CC}{{\rm
C}\kern-5.5pt\rule[.1pt]{.4pt}{1.5ex}\hspace{.5em}} %Complejos

\newtheorem{theorem}{Theorem}

\newtheorem{remark}{Remark}

\begin{document}

\renewcommand\arraystretch{1.0}

\begin{frontmatter}

\title
{Boundary corrections for splitting methods in the time integration of multidimensional parabolic problems}

\author[Lag]{S. Gonz\'{a}lez-Pinto \fnref{aut1}}
\author[Lag]{D. Hern\'{a}ndez-Abreu \fnref{aut1}}

\fntext[aut1]{This work has been partially supported by the Spanish
Project PID2022-141385NB-I00 of Ministerio de Ciencia e Innovaci\'on.}

\cortext[cor1]{Corresponding author: S. Gonz\'{a}lez-Pinto (spinto\symbol{'100}ull.edu.es)}

\baselineskip=0.9\normalbaselineskip
\vspace{-3pt}

%\maketitle

\address[Lag]{ { \footnotesize Departamento de An\'{a}lisis Matem\'{a}tico and
 IMAULL. Universidad de
La Laguna.\\ 38200, La Laguna, Spain 
(spinto\symbol{'100}ull.edu.es, dhabreu\symbol{'100}ull.edu.es).}}

\begin{abstract}
This work considers two boundary correction techniques to mitigate the reduction in the temporal order of convergence in PDE sense (i.e., when both the space and time resolutions tend to zero independently of each other) of $d$ dimension space-discretized parabolic problems on a rectangular domain subject to time dependent boundary conditions. We make use of the MoL approach (method of lines) where the space discretization is made with  central differences of order four and  the time integration is carried out with $s$-stage
AMF-W-methods. The time integrators  are of ADI-type (alternating direction implicit by using a directional splitting) and of  higher order than the usual ones appearing in the literature which only reach order 2. Besides,  the techniques here explained also work for  most of splitting methods, when directional splitting is used. A remarkable fact is that with these techniques, the time integrators recover the temporal order of PDE-convergence  at the level of time-independent boundary conditions. 
\end{abstract}

\begin{keyword}
%% keywords here, in the form: keyword \sep keyword
multidimensional parabolic problem \sep ADI-methods \sep AMF-W method \sep PDE-convergence 
%% MSC codes here, in the form: \MSC code \sep code
%% or \MSC[2008] code \sep code (2000 is the default)
\end{keyword}

\end{frontmatter}

\section{Introduction}
\label{sect:intro}

\subsection{Boundary Corrections  in the time integration of spatially discretized parabolic PDEs}

It has been well known for  sixty years at least that most of one step  methods suffer from order reduction when integrating space discretized (with finite difference, finite volume, finite element, etc)  time-dependent  PDEs in case that  the boundary conditions (BCs), usually Dirichlet, Neumann or Robin BCs, are time dependent. Of course, an initial condition is always required in order to have a well posed PDE problem. For usual one-step methods such as Runge-Kuta or Rosenbrock methods (which  do not make use of splitting), this phenomenon was already studied among others by Sanz-Serna et al.  (1987) \cite{SERNAVER-NUMATH1987} and Carpenter et al.  (1995) \cite{CARPENTER-SISC1995}  by considering explicit Runge-Kutta methods and hyperbolic 1D-PDEs with time dependent Dirichlet BCs. After studying the global errors through the local error propagation and detecting the terms giving rise to order reduction they both propose different remedies to avoid or mitigate such undesirable effect, which involve some modifications of the boundary conditions in  the space discretization at grid-points close to the boundary. For parabolic problems (stiff case) and time integration with implicit Runge-Kutta methods,  the order reduction  for any A($\alpha$)-stable method ($\alpha \ge 0$), decreases in the best situation at $\min\{p,q+1\}$ where $q$ is the stage order of the method and $p$ its classic convergence order. The order reduction phenomenon has also  been analysed for Rosenbrock methods by Rang and Angermann in 2005 \cite{RANG-BIT2005} for DAEs of index 1, by Ostermann and Roche 1993 \cite{OSTROCHE-SINUM1993} by showing  applications in 1D parabolic PDEs mainly and also  by Skvortsov  in several papers from 2014 at present \cite{SVORSTOV-CMMP2014,SVORSTOV-CMMP2022} reaching similar conclusions by using the Prothero-Robinson model \cite{PROTHERO-1974} as a test problem and concluding essentially that the convergence order is limited by $\min\{p,\tilde{q}+1\}$ where $p$ is the classic convergence order of the method and $\tilde{q}$ is the ``stage pseudo-order" of the method which is associated to the global errors associated with the Prothero and Robinson model (recall that the stage order of Rosenbrock methods is limited by one, but its ``pseudo-order" can be arbitrary large \cite{SVORSTOV-CMMP2014,SVORSTOV-CMMP2022}). For  parabolic and quasi-parabolic problems, Lubich and  Ostermann in \cite{LUBOST-MATCOM1995A,LUBOST-APNUM1995B} have also given convergence results which essentially state that in the most of the interior points of the space domain the convergence order of the usual $A(\alpha)$-stable implicit Runge-Kutta methods retains its classical order of convergence, but this order can  decrease  substantially when the grid-points are close enough to the boundary, particularly in case of time-dependent BCs.  The same authors \cite{LUBOST-IMAJNA1995C} have shown that the convergence order is limited by two for  Rosenbrock-type methods and  some class of interesting nonlinear  parabolic problems.

 I. Alonso-Mallo et al. \cite{MALCANO-MATH2021} have considered  a variation of the MoL method where in first place  the time discretization is made by using some Rosenbrock method or  some modifications of it by taking into account the time dependent BCs and some derivatives of them to be included in the internal stages of the method and considering the space discretization after the time discretization. However, with this approach only 1D-parabolic problems  with  a smooth reaction term are presented and their technique does no make use of splitting at all. It allows to recover the classical convergence order of the Rosenbrock method (for some classical  $A(\alpha)$-stable Rosenbrock methods), but only  in the stiff sense, i.e., for prefixed space discretizations based on a fixed number of grid-points, but it is not shown that  the convergence order results hold in PDE-sense, i.e.,  when   the width of the space mesh-grid and  the time step-size tend  to zero  independently of each other. In a more abstract framework,  Calvo and Palencia \cite{CALVOPAL-MATCOM2002} have also treated the order reduction phenomenon for implicit Runge-Kutta methods  applied to parabolic problems  when time-dependent BCs are imposed. They propose some techniques to avoid the order reduction but its actual application to  multidimensional PDE problems seems to be quite difficult and not very cheap due to the fact that none splitting is contemplated. In \cite{CALVO-FRU-NOV-APNUM2001} Linearly implicit Runge-Kutta methods (LIRK methods) of classical orders 3 and 4 are considered for the time integration of advection diffusion reaction PDEs. The authors assume a splitting in two terms where  the diffusion part term (it can also include the advection part in it, but it is assumed a dominant diffusion)  is linear with constant coefficients and the other term could be nonlinear but smooth and it includes the reaction part and/or the advection. In spite of  the fact that the results presented show that the methods are  competitive  with  some implementation of the BDF method (Backward Differentiation Formulae) in two parabolic problems (1D and 2D respectively) with homogeneous BCs, it is not clear (and not demonstrated) that the PDE convergence order coincides with the classical convergence order of the methods and that the results are competitive for the case of time-dependent BCs. 

For splitting methods of one-step nature that make use of directional splitting in  multidimensional $dD$-parabolic PDEs ($2D,\: 3D,\:...$), by taking into account their weaker  stability properties than  those of the implicit  Runge-Kutta or Rosenbrock families,  it is expected more dramatic reductions in their convergence orders than in the latter ones, which  have the advantage of not  making use of splitting, but its practical implementation is almost prohibitive in $dD$ parabolic problems for $d>2$.  However, the splitting methods are quite attractive due to its low computational costs compared with Runge-Kutta, Rosenbrock or Linear Multistep methods due to the fact that the former ones allow to carry the algebra costs  at the level of 1D-problems multiplied by the number of spatial dimensions $d$ and  at the same time to have  good enough stability properties to provide long time integrations with reasonable large time step-sizes when  appropriate  methods of splitting-type are chosen \cite[Chapt. IV]{HV-2003}, \cite{GHH-JCAM2023,GHH-ESAIM2023,HOWSOM-JCAM2001}. 

  Next,  we just mention a few relevant works where it is  observed the effect of the order reduction for splitting methods (parabolic problems) and  some remedies  to avoid or mitigate it. In all cases, as far as we know,  order  of convergence two in the time integration in PDE-sense (that is, global errors are bounded in some $\ell_p$-norm by $C_1 (\Delta t)^2 + C_2h^2$, where $C_1$ and  $C_2$ are two constants independently of  $h\rightarrow 0$ and $\Delta t\rightarrow 0$,  where $h$ denotes the width of the spatial mesh-grid and $\Delta t$ the time step-size. Here, it is assumed a second order spatial discretization based in the standard central differences) seems to be a barrier for the most interesting splitting methods proposed in the literature \cite[Chap. IV]{HV-2003}, \cite{MARCHUK-1990}, even in the most favourable case of homogeneous BCs.   We next mention some interesting papers about splitting methods and order reduction by time-dependent boundary conditions. Fairwheater and Mitchel  \cite{FAIRMIT-SINUM1967} treat  ADI-methods  (Alternating Direction Implicit) for the Douglas scheme \cite{DOUGLAS-1955,DOUGLASR-1956} in 2D-problems by considering directional splitting. In Sommeijer et al.  \cite{SOMHOUWVER-IJNME1981,HOWSOM-JCAM2001} some  LOD-methods (Locally One Direction Methods) and ADI-methods with directional splitting in  2D and 3D problems are dealt with. Other remarkable works are those ones by  LeVeque (1985) \cite{LeVEQUE-1985} by considering  the LOD-CN method (Locally One Direction Crank-Nicolson method also known as Yanenko's method) with directional splitting, and  Marchuck  \cite{MARCHUK-1971} also considers and collects results about several  ADI methods and LOD methods with directional splitting, but without proving theoretical results about the PDE-convergece order of the methods. A very nice survey about the order reduction phenomenon until the beginnings of the 21st century can be seen in the excellent monograph by Hundsdorfer and Verwer (2003) \cite[Chapt. IV]{HV-2003} and some extense treatment  about methods based on AMF-splitting  (Approximate Matrix Factorization) is \cite{HOWSOM-JCAM2001}.
 
 Another famous kind of splitting methods are those ones so-called IMEX methods (Implicit-Explicit methods), thus  among many others, we mention that Asher et al. at the end of the nineties  \cite{ASHERRS-1997,ASHERRS-1995},  considered one step IMEX methods (Implicit-Explicit splitting)  by using a two term splitting, one for the stiff part (implicit treatment) and another one for the explicit part, but they have neither considered the effect of the order reduction due to time dependent BCs nor  the PDE-convergence order. Crouzeix in 1980 \cite{CROUZEIX-NUMATH1980} and Arrarás et al. in 2021 \cite{HOUTHUNPOR-BIT2017} also analyzed the linear multistep  approach of IMEX type for parabolic problems, but as we know multistep methods suffer from order reduction due to its limited stability properties which are more pronounced in the case of using splitting (do not suffer order reduction due to the consistency order). Additionally, Hundsdorfer et al. have made significant contributions in the stability analysis and convergence study of the Douglas  method, the ADI-methods and  a particular splitting of the  Trapezoidal Rule  as it can be seen e.g. in \cite{HV-2003,HV-MATCOM1989,HUNDS-MATCOMP1992,HUNDS-APNUM2002}. However,  for all these one step methods the time PDE-convergence order was limited to two. 
 
 Other kind of splitting, often call operator splitting, consists of splitting the original PDE problem in  two (or more) PDE subproblems with appropriate boundary conditions in such a way that their solutions  adequately combined  provide a consistent  approximation to the original PDE problem. Particularly famous are the Strang Splitting (order two) and the Lie-Trotter splitting (order one). In several papers, Alonso-Mallo et al. \cite{MALCANOREG-2021} and Einkemmer et al.   \cite{EINKOST-SISC2015,EINKOST-SISC2016,EINMOCOST-AMC2018}, by considering the Strang splitting  on diffusion reaction PDE problems with time-dependent BCs, appreciated an order reduction until order close to one,  in most of the interesting $\ell_p$-norms, but not in the case of homogeneous BCs where the order two was preserved. To overcome the order reduction in case of time-dependent BCs, the authors used an splitting in two terms,  the diffusion part and  the reaction part with appropriate boundary conditions for each term. Einkemmer et al. \cite{EINKOST-SISC2015,EINKOST-SISC2016,EINMOCOST-AMC2018} called the new method the modified Strang scheme. This method recovers  the original order of convergence two in the most interesting $\ell_p$-norms, namely $\ell_1,\:\ell_2$ and maximum norm and for different types of time-dependent boundary conditions. They also considered some modification of the Lie-Trotter scheme which reached order two in a very few particular situations. To our point of view, despite  its value for some kind of interesting problems, it cannot be widely used due to the limitation that the diffusion term has to be treated as a whole and  no splitting  is contemplated in it. Additionally the  time PDE convergence order is limited to two.  
 
A. Arrar\'as and L. Portero in \cite{ARRARAS-POR-AMC-2015} have considered a  domain-decomposition splitting  consisting in splitting  the whole space  domain in small sub-domains and to apply some modified Crank-Nicolson to each subdomain  by taking into account the transference of information on the boundaries of the overlapping subdomains and adding some corrector terms of order $\mathcal{O}(\Delta t^3)$. The algorithm is based on the Douglas method, it is well described and some numerical illustrations on some  2D problems seem to confirm PDE convergence of order two, but as the authors there recognize the convergence of order two  was far to be proved at that time.
 
 From the reasons above exposed it seems to be clear that directional splitting in the basis of a MoL approach is an attractive alternative to deal with multidimensional problems  (it has been widely used in the literature, see e.g.  \cite{HV-2003} and references therein) since it is quite  straightforward to apply, but has the drawback that the case of  time-dependent boundary conditions involves order reduction  in the PDE convergence orders of the methods. The goal of this paper is to provide some techniques to avoid or mitigate the order reduction at the level of time-independent  boundary conditions for the troublesome case of time-dependent boundary conditions and also to show that PDE order of convergence three in time is possible even in the maximum norm, which is the most problematic case of the $\ell_p$ norms. The techniques here described are designed for AMF-W methods but they can be applied to most of standard splitting methods with directional splitting  to get  similar positive effects.  AMF-W methods with directional splitting are of 1D-nature and allow obtain PDE order of convergence  three for  $dD$ parabolic problems  with $d $ arbitrary \cite{GHH-JCAM2023,GHH-ESAIM2023}.    We are not aware of other methods of splitting type having such higher PDE orders of convergence in time.
 
 The remainder of the paper is organized as follows, taking into account that we follow the MoL approach. In Section 2 we introduce the space discretization (and the splitting to be used) based on finite differences, which turns out to be fourth order and not second order as it is usual. We recall the format of the AMF-W methods (for the time integration) and we complete the section by illustrating the  order reduction phenomenon (in two norms, the weighted $\ell_2$-norm and the maximum norm) on an academic 3D parabolic PDE problem. In Section 3 we include two convergence theorems that explain most of the observed numerical results of the previous section. In section 4 we provide a technique to avoid the convergence order reduction by using some adequate ``interpolant" that fulfils the time dependent boundary conditions. In section 5 we propose another practical technique to mitigate the order reduction which is based on some extension of the elliptic operator of the PDE to the boundary and could be more easily implemented than the case of the ``interpolant". In section 6, the effectiveness of the techniques are illustrated with several 2D and 3D linear and nonlinear problems. Finally, some conclusions and future research are drawn in section 7. It is also worth to mention that most of the results here presented have been developed along  the last seven years by Gonzalez-Pinto and collaborators  in  \cite{GHHP-SISC2018,GHH-SINUM2020,GHPSRS-JCPHYS2022,GHP-JCAM2021,GHP-NUMAL2023}.

 \section{The PDE problems considered}\label{sect-2}

We consider semilinear parabolic PDEs ($d$ space dimensions) with a possible nonlinear reaction term $r$ (a nonstiff or mildly stiff reaction), where below  $\mathcal{L}$ represents an   ellyptic differential operator:

\begin{equation}\label{eq1-0}\begin{array}{ll}
\hbox{\rm (PDE)} & 
 u_t(\vec{x},t)=\mathcal{L} u + r(\vec{x},t,u), \quad  \vec{x} \in \Omega=(0,1)^d, \quad t \in I^*=(0,t^*], \\[0.4pc] &
\mathcal{L} u(\vec{x},t)=\sum_{j=1}^d \Big( a_j(\vec{x},t) \: \partial_{x_jx_j} u + b_j(\vec{x},t)  \: \partial_{x_j} u\Big) \\[0.4pc] & 
a_j(\vec{x},t)\ge \bar{a}>0, \; (\vec{x},t)\in \Omega\times I^*, \;j=1,\ldots,d. \\[0.4pc] \hbox{\rm (BCs)} & 
 \gamma(\vec{x},t,u,\partial_{x_1}u,\ldots,\partial_{x_d}u)=\beta(\vec{x},t), \quad  
\vec{x} \in \partial \Omega, \quad t \in I^*,\\[0.4pc]  & \gamma \; \hbox{\rm linear on  } u,\:\{\partial_{x_j} u\}_{j=1}^d. \; \hbox{\rm Typically Dirichlet, Neumann or Robin BCs.}\\[0.4pc]  
(IC) & u(\vec{x},0)=u_0(\vec{x}), \quad \vec{x} \in  \Omega. \end{array}
\end{equation}
 Of course, it is assumed regularity of the data up to the boundary (including it), i.e., the coefficients of the linear operator, the reaction term, the initial condition, the  boundary conditions and the exact solution $u$  have continuous derivatives up to a certain order on $\bar{\Omega}\times I^*$, with $\bar{\Omega}=\Omega\cup \partial{\Omega}$ and $\partial{\Omega}$ denoting the boundary of $\Omega$.
 
To fully discretize the PDE problem, we follow the MoL approach (method of lines), where the space discretization is carried out by using finite differences on evenly spaced nodes  (the conservative form of the operator, namely  $\mathcal{L}u=\sum_{j=1}^d \partial_{x_j}\Big( a_j(\vec{x},t) \: \partial_{x_j} u +  b_j(\vec{x},t)  u\Big)$, could also be used  combined with finite differences, but for simplicity we omit that in this paper). For the time integration we apply  AMF-W-methods, which will make use of a directional splitting as it  will be described below. 

%%%%%%%%%%%%%%%%%%%%%%%%%%%%%%%%%%%
%%%%% 

 \subsection{The Spatial Semidiscretization}\label{sect1-4}

Space discretizations based on central differences yield differential systems of ODEs of the form 
\begin{equation}\label{eq1-1} \begin{array}{l} \hbox{\rm (ODE) }\quad \displaystyle \dot V(\vec{x}_G,t)=\mathcal{L}^{(h)} V + r^{(h)}(\vec{x}_G,t,V) ,\quad  V(\vec{x}_G,t)\approx u(\vec{x}_G,t)\\[0.4pc]
\hbox{\rm (BCs)}\;\quad \gamma^{(h)}(\vec{x}_G,t,V)=\beta(\vec{x}_G,t), \quad  
\vec{x} \in \partial \Omega_h, \quad t \in I^*=(0,t^*]\\[0.4pc] 
\hbox{\rm (IC) } \qquad V(\vec{x}_G,0)=u_0(\vec{x}_G), \quad \vec{x}_G \in  \Omega_h.\end{array}\end{equation}
Here $\gamma^{(h)}$ is the discretized boundary operator which is linear on V, where typically Dirichlet, Neumann or Robin BCs are imposed. We  consider the case of evenly spaced nodes and make use of the following notations,
$$  \begin{array}{l}\vec{x}_G =(x_1^{(j_1)}, x_2^{(j_2)},\ldots,x_d^{(j_d)})^\top \in \Omega_h, \quad t \in I^*,\quad \hbox{\rm where }\\[0.4pc]
 \displaystyle x_l^{(j_l)}=j_l \Delta x_l, \quad \Delta x_l=\frac{1}{1+n_l},\quad j_l=0,1,2,\ldots,n_l+1;\quad l=1,2,\ldots,d.\end{array}$$
The space discretization of the linear operator $ \mathcal{L}$ is indicated in (\ref{eq1-4})  where the operators for the diffusion and advection terms are based on central differences and will be defined below,
 \begin{equation}\label{eq1-4}
 \mathcal{L}^{(h)} V(\vec{x}_G,t)=\sum_{j=1}^d \Big( a_j(\vec{x}_G,t) \: \partial^{(h)}_{x_jx_j} V + b_j(\vec{x}_G,t)  \: \partial^{(h)}_{x_j} V\Big).
 \end{equation}
We recall that $\Omega_h$ will denote the interior grid-points (grid-points in $\Omega$)  for Dirichlet boundary conditions. For other boundary conditions such as Neumann or Robin conditions in the whole or part of the boundary,  $\Omega_h$ would include the corresponding  boundary points too. 
 
\begin{figure}
\begin{center}

\begin{tikzpicture}[x=1cm,y=1cm,scale=0.6]

\draw[thick,red] (0,0)--(2,0);
\draw[thick] (2,0)--(3,0);
\draw[dashed] (3,0)--(5,0);
\draw[thick] (5,0)--(11,0);
\draw[dashed] (11,0)--(13,0);
\draw[thick] (13,0)--(14,0);
\draw[thick,red] (14,0)--(16,0);

\draw[thick,red] (0,-0.2)--(0,0.2);
\draw[thick,red] (1,-0.2)--(1,0.2);
\draw[thick,red] (2,-0.2)--(2,0.2);
\draw[thick] (6,-0.2)--(6,0.2);
\draw[thick] (7,-0.2)--(7,0.2);
\draw[thick] (8,-0.2)--(8,0.2);
\draw[thick] (9,-0.2)--(9,0.2);
\draw[thick] (10,-0.2)--(10,0.2);
\draw[thick,red] (14,-0.2)--(14,0.2);
\draw[thick,red] (15,-0.2)--(15,0.2);
\draw[thick,red] (16,-0.2)--(16,0.2);
\draw[thick,red] (1,0) node[]  {\normalsize x};
\draw[thick] (8,0) node[]  {\normalsize x};
\draw[thick,red] (15,0) node[]  {\normalsize x};
\draw [] (0,-0.2) node[below]  {\scriptsize $0$};
\draw [] (1,-0.2) node[below]  {\scriptsize $x_1$};
\draw [] (2,-0.2) node[below]  {\scriptsize $x_2$};
\draw [] (6,-0.2) node[below]  {\scriptsize $x_{i-2}$};
\draw [] (7,-0.2) node[below]  {\scriptsize $x_{i-1}$};
\draw [] (8,-0.2) node[below]  {\scriptsize $x_i$};
\draw [] (9,-0.2) node[below]  {\scriptsize $x_{i+1}$};
\draw [] (10,-0.2) node[below]  {\scriptsize $x_{i+2}$};
\draw [] (14,-0.2) node[below]  {\scriptsize $x_{M-1}$};
\draw [] (15,-0.2) node[below]  {\scriptsize $x_M$};
\draw [] (16,-0.2) node[below]  {\scriptsize $1$};

\draw[thick,red] (0,-3)--(2,-3);
\draw[thick,red] (0,-3.2)--(0,-2.8);
\draw[thick,red] (1,-3.2)--(1,-2.8);
\draw[thick,red] (2,-3.2)--(2,-2.8);

\draw [red] (-1,-3.7) node[]  {\small $\frac{1}{\Delta x^2} \cdot $};

\draw [red] (0,-3.2) node[below]  {\scriptsize $1$};
\draw [red] (1,-3.2) node[below]  {\scriptsize $-2$};
\draw [red] (2,-3.2) node[below]  {\scriptsize $1$};

\draw[thick] (6,-3)--(10,-3);

\draw[thick] (6,-3.2)--(6,-2.8);
\draw[thick] (7,-3.2)--(7,-2.8);
\draw[thick] (8,-3.2)--(8,-2.8);
\draw[thick] (9,-3.2)--(9,-2.8);
\draw[thick] (10,-3.2)--(10,-2.8);

\draw [] (4.8,-3.7) node[]  {\small $\frac{1}{\Delta x^2} \cdot $};

\draw [] (6,-3.2) node[below]  {\scriptsize $\frac{-1}{12}$};
\draw [] (7,-3.2) node[below]  {\scriptsize $\frac{16}{12}$};
\draw [] (8,-3.2) node[below]  {\scriptsize $\frac{-30}{12}$};
\draw [] (9,-3.2) node[below]  {\scriptsize $\frac{16}{12}$};
\draw [] (10,-3.2) node[below]  {\scriptsize $\frac{-1}{12}$};

\draw[thick,red] (14,-3)--(16,-3);
\draw[thick,red] (14,-3.2)--(14,-2.8);
\draw[thick,red] (15,-3.2)--(15,-2.8);
\draw[thick,red] (16,-3.2)--(16,-2.8);

\draw [red] (13,-3.7) node[]  {\small $\frac{1}{\Delta x^2} \cdot $};

\draw [red] (14,-3.2) node[below]  {\scriptsize $1$};
\draw [red] (15,-3.2) node[below]  {\scriptsize $-2$};
\draw [red] (16,-3.2) node[below]  {\scriptsize $1$};

\end{tikzpicture}%
\caption{Stencils used on each direction for the diffusion case (operator $\partial^{(h)}_{x_j,x_j}V$). Observe that they are based on central second order differences for the adjacent points to the boundary and in  central fourth order differences for the other interior points. For the advection terms (operator $\partial^{(h)}_{x_j}V$) the corresponding approaches of second and fourth order will be taken at the same points.}\label{Fig-1}
\end{center} 
\end{figure}
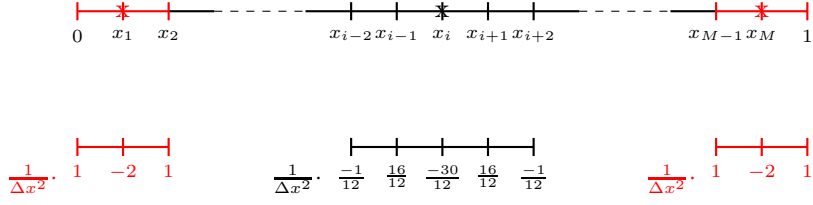

We will employ fourth order space discretizations based on fourth order central differences for the interior points  no-adjacent to the boundary and a second order discretization for the interior points adjacent to the boundary (in case of Dirichlet BCs), i.e., on each direction (spatial variable) we take the stencils given in Figure \ref{Fig-1} for the diffusion (the corresponding ones should be taken for the advection terms), which still gives a general fourth order approach in space in the weighted $\ell_2$-norm as we will see in the numerical examples later considered in section \ref{sect-5}.  Of course, we can choose second order central differences as it is usual, but due to the fact that we are considering  high order methods  for the time integration (order three or more), it is reasonable to use high order discretizations for the space too. 

%%%%%%%
In the case of Dirichlet boundary condition the problem (\ref{eq1-1}) reduces to 
 \begin{equation}\label{eq1-20} \displaystyle \dot V(t)=F(t,V)=\sum_{j=0}^d F_j(t,V),\quad V(0)=V_0, \quad V(t)\equiv V(\vec{x}_G,t)\end{equation}
where we consider the case in which the splitting is directional with the reaction term ($F_0$) set in the first place. This is adventageous to mitigate the reduction in the convergence order when using AMF-W-methods as it was indicated in   \cite{GHH-JCAM2023,GHH-ESAIM2023}, 
\begin{equation}\label{eq1-21}\begin{array}{l}
 F_0(t,V)= r^{(h)}(\vec{x}_G,t,V),\\[0.5pc]
 F_j(t,V)= a_j(\vec{x}_G,t) \: \partial^{(h)}_{x_jx_j} V + b_j(\vec{x}_G,t)  \: \partial^{(h)}_{x_j} V, \qquad (j=1,2,\ldots,d).
\end{array} \end{equation}
It should be observed that the dimensions of the terms of the ODE problem are
$$\mbox{\rm dim } V=\mbox{\rm dim } F_j=\mbox{\rm dim } F= n_1\times n_2\times \ldots \times n_d=N_G, \;(j=1,2,\ldots,d).$$

\medskip

\subsection{$s$-stage AMF-W method $(A,L,b,\theta)$} 

We recall that a Rosenbrock method $(A,L,b,\theta)$ for the integration of ODE systems (\ref{eq1-20}) can be rewritten as (it is not a splitting method), see e.g. \cite{HW-2},
\begin{equation} \label{eq1-7a} 
\begin{array}{rcl}
(I-\theta \Delta t W_n) K_i&=&\Delta t F(t_n+c_i\Delta t, V_n+\sum_{j=1}^{i-1} a_{ij} K_j)+ \sum_{j=1}^{i-1} l_{ij} K_j \\[0.5pc]&+ & \theta \rho_i \Delta t^2 \dot F(t_n,V_n)\quad i=1,2,\ldots,s\\[0.5pc]
V_{n+1}&=&V_n+\sum_{i=1}^s b_i K_i,
\end{array}
\end{equation} 
where 
\begin{equation}\label{eq1-8} (c_i)_{i=1}^s=A (\rho_i)_{i=1}^s,\quad  (\rho_i)_{i=1}^s=(I_s-L)^{-1} \hbox{\bf 1}, \quad \hbox{\bf 1}=(1,1,\ldots,1)^\top \in \mathbb{R}^s,
\end{equation}
$\dot F(t,V)=\partial_t F(t,V) $ and $W_n=\partial_V F(t_n,V_n)$ or  $W_n=\partial_V F(t_n,V_n)+\mathcal{O}(\Delta t).$ 
In the latter case, the method is sometimes called a W-method \cite{STEIHAUG-WOLFB-1979} of particular type.

Due to the high computational costs involved with the integration with Rosenbrock methods for problems (\ref{eq1-20}) stemming from $dD$ multidimensional PDEs ($d\ge 2$), we are going to consider  $s$-stage  AMF-W methods $(A,L,b,\theta)$, which are a cheaper simplification of Rosenbrock methods and  present the following formulation for a directional splitting based in $d+1$ terms \cite{GHHP-SISC2018} 
\begin{equation}\label{eq1-7}
\begin{array}{rl}\hbox{\rm for } i=1,2,\ldots,s, & \hbox{\rm compute the stages from}\\ 
K_i^{(-1)}=&\hspace*{-2mm}\Delta t \:F\big(t_n+c_i\Delta t,V_n+\displaystyle{\sum_{j=1}^{i-1} a_{ij} {K}_j\big) + \sum_{j=1}^{i-1} \ell_{ij}}{K}_j, \\[2mm] 
(I-\theta \Delta t D_j)K_i^{(j)}=& K_i^{(j-1)}+  \theta \rho_i\Delta t^2 \dot 
F_j(t_n,V_n), \quad (j=0,1,\ldots,d)\\[2mm]
{K}_i=&{K}^{(d)}_i, \qquad D_j=\partial_V F_j(t_n,V_n), \quad \dot 
F_j=\partial_t F_j(t_n,V_n).\end{array}
\end{equation}
The step-point solution $V_{n+1} \approx u(\vec{x}_G,t_{n+1})$ is computed from  
\begin{equation}\label{eq1-8a} V_{n+1}=V_n + \sum_{i=1}^s b_i K_i.\end{equation}

\subsection{Illustration of the reduction in the temporal PDE-convergence order on a  3D parabolic academic example}
As a test example we use the following 3D parabolic PDE with Dirichlet BCs,  a source term $r$ and an initial condition so that the exact solution is given by 
\begin{equation}\label{eq1-10}\displaystyle u=\beta(\vec{x},t)= 4^3e^t \prod_{j=1}^3 x_j(1-x_j)+ C\:e^t \sum_{j=1}^3 \Big(x_j+\frac{1}{j+2}\Big)^2. \end{equation}
We then consider, 
\begin{equation}\label{eq1-11} \begin{array}{l}
\hbox{\rm (PDE): } \displaystyle u_t(\vec{x},t)=\sum_{j=1}^3  \partial_{x_jx_j} u + r(\vec{x},t), \quad (\vec{x},t) \in \Omega\times I^*, \quad  \Omega=(0,1)^3,\\[0.5pc]
\hbox{\rm (BCs): } u(\vec{x},t)=\beta(\vec{x},t), \quad  
(\vec{x},t) \in \partial \Omega\times I^*, \quad I^*=(0,1],\\[0.5pc]
\hbox{\rm (IC): } \; u(\vec{x},0)=\beta(\vec{x},0), \quad \vec{x} \in  \Omega. \end{array} \end{equation}
The source term is given by,  
\begin{equation}\label{eq1-12} r(\vec{x},t):=\dot \beta(\vec{x},t)-\sum_{j=1}^3 \partial_{x_jx_j} \beta(\vec{x},t). \end{equation}
It should be noted that in case that $C=0$, we have homogeneous Dirichlet BCs and in case of $C=1$, time dependent Dirichlet BCs are imposed. The problem is discretized in space by using second order central differences at adjacent points to the boundary and by  fourth order central differences at  the other interior points, as it was previously indicated. Observe that  no space errors are present in the semidiscretized ODE systems due to the polynomial form of the exact solution on the space variables.

\begin{table}
\begin{center}
\begin{tabular}{l|l|l} $\Delta t =\Delta x_j$ & \ GE$_{2,h} \quad$ ($p_2$) & \ GE$_\infty \qquad$ ($p_\infty$) \\ \hline 
$1/4 $ & 0.33e+0 \;(---) &  0.52e+0 $\quad$(---) \\
$1/8 $ & 0.60e-1 \;(2.45) &  0.11e+0 $\;$ (2.19) \\
$1/16 $ & 0.97e-2 \;(2.64) &  0.20e$\:$-1 $\;\:$ (2.64) \\$1/32 $ & 0.14e-2 \;(2.74) &  0.30e$\:$-2 $\;\:$ (2.75) \\$1/64 $ & 0.20e-3 \;(2.84) &  0.39e$\:$-3 $\;\:$ (2.91) \\$1/128 $ & 0.27e-4 \;(2.90) &  0.49e$\:$-4 $\;\:$ (3.01) \\ \hline
\end{tabular} \caption{Case $C=0$ in the test problem (\ref{eq1-10})-(\ref{eq1-12}). Global errors at the end-point and observed convergence orders (in parenthesis) for the corresponding norms.} 
\label{table1}
\end{center}
\end{table}

As time integrator we consider the AMF-W-method proposed in \cite[p. 400]{HV-2003}
\vspace{0.5cm}
\begin{equation} \label{eq1-9} \begin{array}{c}
 A=\left(\begin{array}{cc} 0 & 0 \\ 2/3 & 0 \end{array} \right), \quad  L=\left(\begin{array}{cc} 0 & 0 \\ -4/3 & 0 \end{array} \right),\quad b= \left(\begin{array}{c} 5/4 \\ 3/4  \end{array} \right) 
   \\[0.7pc] 
 \displaystyle{ \theta =\frac{3+\sqrt{3}}{6}.}\end{array} \end{equation}
The method has order of convergence three in time in ODE sense  for a fixed space resolution (fixed ODE), but this is not guaranteed in the case that both $$\Delta x_j \rightarrow 0, \quad \Delta t \rightarrow 0 \qquad \hbox{\rm (PDE order of convergence).} $$
We  numerically check the observed orders of convergence (PDE convergence orders) on the 3D academic example above for the particular case in that   $$h=\Delta x_j=\Delta t\rightarrow 0, \quad j=1,2,3.$$
By using the fact that no-space discretizations errors of the PDE arise in the ODEs,
the global errors $GE(t^*,h,\Delta t):=u(\vec{x}_G,t^*)-V(\vec{x}_G,t^*) $, are bounded in the weighted $\ell_2$-norm  \cite{GHH-ESAIM2023}
\begin{equation}\label{eq1-15} \displaystyle  \| V \|_{2,h}^2= \Delta x_1 \Delta x_2 \cdots \Delta x_d \sum_{j_l=1\atop l=1,2,\ldots,d }^{n_l} |V(x_1^{(j_1)},x_2^{(j_2)},\ldots,x_d^{(j_d)})|^2,\end{equation}
as 
\begin{equation}\label{eq1-13a}
\|GE(t^*,h,\Delta t)\|_{2,h} \le C^* (\Delta t)^p, \quad \forall h>0, \Delta t>0,
\end{equation}
where $C^*>0$ is a constant independently of $\Delta t$ and $\Delta x_j$ for $j=1,\ldots,d,$ and $p>0$ is the temporal order (in PDE sense) of the AMF-W method. 
Additionally to the bound above, it is expected that the global errors of any convergent AMF-W method for $\Delta t= h$  behave as 
\begin{equation}\label{eq1-13b}
GE(t^*,h,h)= \chi_V(\vec{x}_G,t^*) h^p(1+ \mathcal{O}(h)^{r}), \quad p>0, \;r>0,
\end{equation}
with $\chi_V(\vec{x}_G,t^*)$ being a vector with the grid-values of some smooth function $\chi(\vec{x},t^*)$ defined on $\bar \Omega\times I^*$. In such a situation, for the weighted $\ell_2$-norm (\ref{eq1-15}) or the maximum norm,   it is expected that 
\begin{equation}\label{eq1-14}
\|\chi(\vec{x}_G,t^*)\| \simeq \kappa,\; \hbox{\rm a constant, independently of the space mesh-grid } \Omega_h.
\end{equation}
From this fact, for each $h=\Delta t=\Delta x_j,\;j=1,2,3$, we have estimated the PDE-convergence orders in  Table \ref{table1} and Table \ref{table2} in the following way 
\begin{equation}\label{eq1-16}
p\simeq \log (\|GE(t^*,2h,2h)\|/ \|GE(t^*,h,h)\|) /\log 2.   
\end{equation}

  Observe that in case of homogeneous Dirichlet BCs ($C=0$) no order reductions in the convergence orders are observed either for the weighted $\ell_2$-norm   or in the maximum norm 
$$\displaystyle  \| V \|_{\infty}= \max_{1\le j_l\le n_l \atop l=1,2,3} |V(x_1^{(j_1)},x_2^{(j_2)},x_3^{(j_3)})|.$$

\begin{table}
\begin{center}
\begin{tabular}{l|l|l} $\Delta t=\Delta x_j$ & \ GE$_{2,h} \quad$ ($p_2$) & \ GE$_\infty \qquad$ ($p_\infty$) \\ \hline 
$1/4 $ & 0.40e+0 \;(---) &  0.96e+0 $\quad$(---) \\
$1/8 $ & 0.70e-1 \;(2.53) &  0.17e+0 $\;$ (1.68) \\
$1/16 $ & 0.12e-1 \;(2.55) &  0.98e$\:$-1 $\;\:$ (0.80) \\$1/32 $ & 0.22e-2 \;(2.41) &  0.59e$\:$-1 $\;\:$ (0.74) \\$1/64 $ & 0.48e-3 \;(2.23) &  0.34e$\:$-1 $\;\:$ (0.78) \\$1/128 $ & 0.11e-3 \;(2.11) &  0.20e$\:$-1 $\;\:$ (0.81) \\ \hline
\end{tabular}
\caption{Case $C=1$ in the test problem (\ref{eq1-10})-(\ref{eq1-12}). Global errors at the end-point and observed convergence orders (in parenthesis) for the corresponding norms.} 
\label{table2}
\end{center}
\end{table}

However, in case of time dependent Dirichlet BCs ($C=1$), from Table \ref{table2} it can be appreciated  that 
\begin{enumerate}
\item in the weighted $\ell_{2,h}$-norm,  a little loss in accuracy and  an order reduction in about one unit in the convergence is observed. 
\item In the $\ell_\infty$-norm a big loss in accuracy and a dramatic order reduction (2 units at least) can be seen. 
\end{enumerate} 
Some theoretical justifications about the PDE-convergence orders observed in Tables \ref{table1}-\ref{table2} will be given in the next section. 

\section{Some theoretical convergence results}
%\colorlet{shadecolor}{blue!15} \vspace*{-3mm}
%\begin{shaded}
The next two theorems apply to the case of Dirichlet Boundary conditions for the weighted $\ell_2$ and the maximum norms respectively and they were adapted from those ones given in \cite{GHH-ESAIM2023,GHH-SINUM2020}. They provide some justification for the numerical findings in Table \ref{table1} and Table \ref{table2}, respectively. The first theorem below applies to the weighted Euclidean norm and the second one to the Maximum norm. 

We recall that the stability function of a AMF-W-method $(A,L,b,\theta)$ is obtained from applying the method to the Dahlquist test problem with splitting  \cite[Chapt. IV]{HV-2003}
\begin{equation} \label{eq2-1}
y'(t)=\lambda y, \quad y(t_n)=y_n,\quad \lambda=\sum_{j=0}^d \lambda_j, \quad \lambda_j \in \mathbb{C}, \quad z_j=\Delta t \:\lambda_j.
\end{equation}
Here, the splitting is made in the following way, $F_0(t,y)=\lambda_0y=0, \;F_j(t,y)=\lambda_j y$.
The method gives the  following relation $y_{n+1}=R(z_1,z_2,\ldots,z_d) y_n,$ where the stability function (below, $I_s$ denotes the identity matrix of dimension $s$)  
\begin{equation} \label{eq2-1a}
R(z_1,z_2,\ldots,z_d)= 1+ z b^\top (\Pi \:I_s-L-zA)^{-1}\hbox{\bf 1}, \quad z=\sum_{j=1}^d z_j,\;\Pi=\prod_{j=1}^d(1-\theta z_j), \end{equation}
which  is relevant for linear PDEs with constant coefficients discretized with  standard central differences \cite[Chapt. IV]{HV-2003}, 
\begin{equation} \label{eq2-2} \dot u(t)=\sum_{j=1}^d (J_j u + g_j(t)), \quad t \in [0,t^*].\end{equation} 

  \begin{theorem}\label{th-1} {\rm (Weighted Euclidean norm, \cite{GHH-ESAIM2023})}. Consider  linear systems of  type (\ref{eq2-2}) obtained from the space discretization of constant coefficient PDEs of type (\ref{eq1-0}) with $a_j>0,\;b_j=0,\;j=1,\ldots,d$, by using second order central differences and  Dirichlet boundary conditions. Additionally, assume  that 
\begin{enumerate}
\item    $u(\vec{x},t)$ admits bounded partial  derivatives up to order $p+1$ in  $ \bar 
\Omega \times I^*$. \item  The AMF-W-method $(A,L,b,\theta)$ has order $p$ as  Rosenbrock method, see (\ref{eq1-7a}), whenever \   $W_n-\partial_V F(t_n,V_n)=\mathcal{O}(\Delta t)$ in ODE sense, i.e. negative powers of  $h=\max_{1\le j \le d}\Delta x_j$, can be hidden in the $\mathcal{O}(\Delta t)$ term.
\item The linear stability function $R(z_1,z_2,\ldots,z_d)$ in (\ref{eq2-1a}), satisfies\vspace{-0.3cm}  $$ \begin{array}{c}
 -1 \le R(x_1,x_2,\ldots,x_d) \le 1+ \displaystyle C\frac{x_1+\ldots + x_d}{(1-\theta x_1)\ldots(1-\theta x_d)},  \\[0.5pc]  \quad \forall x_j<0 \;(j=1,2,\ldots,d), \quad C>0 \; \hbox{\rm is a constant.} 
\end{array}$$
\end{enumerate}  
Then, for the global errors $\epsilon_h(t_n):=u(\vec{x}_G,t_n)-V_h(\vec{x}_G,t_n)$, there exist two positive constants $C_1$ and $C_2$ such that, \\ {\rm (a) }    
$$ \Vert \epsilon_h(t) \Vert_{2,h}\le C_1 \Delta t^{\min \{p,3.25^*\}} +C_2 h^2; \quad (\Delta t\rightarrow 0, h\rightarrow 0),$$
in case of time-independent BCs\footnote{Henceforth, PDE order $q^*$ means that any PDE order smaller than $q$ can be reached, but not exactly PDE order $q$.} \\
{\rm (b) }
$$\begin{array}{c}\displaystyle \Vert \epsilon_h(t) \Vert_{2,h}\le C_1 \Delta t^{\min\{p,2\}}+C_2 h^2;\quad (\Delta t\rightarrow 0, h\rightarrow 0), 
 \end{array}  
$$
in case of time-dependent BCs.
 \end{theorem} 
{\sf Proof.} It follows from Theorems 3.1, 3.2 and 3.3 in \cite{GHH-ESAIM2023}. \hfill $\Box$

\subsection{Convergence results in the $\ell_\infty$-norm}
%\colorlet{shadecolor}{blue!15} \vspace*{-3mm}
%\begin{shaded}
 \begin{theorem}\label{th-2} {\rm (Maximum norm, \cite{GHH-SINUM2020,GH-APNUM2022})}. Under the same space discretizations, assumptions 1) and 2) of Theorem \ref{th-1} and  
the stability condition (we take as above $h=\max_{1\le j \le d}\Delta x_j$) 
$$\|R(\Delta t J_1, \ldots, \Delta t J_d)^n \|_\infty \le C, \quad \Delta t \rightarrow 0, \; h\rightarrow 0, \;n=1,2\ldots,t^*/\Delta t,$$ Then,  there exist two positive constants $C_1$ and $C_2$ such that, \\ 
{\rm (a) }  
$$ \Vert \epsilon_h(t) \Vert_{\infty}\le C_1 \Delta t^{\min \{p,2\}} +C_2 h^2; \quad (\Delta t\rightarrow 0, h\rightarrow 0),$$
 in case of time-independent BCs in case of time-dependent BCs\footnote{This result was only shown for A-stable one-stage AMF-W methods of ODE order 2}.\\
{\rm (b) } 
$$\begin{array}{c}\displaystyle \Vert \epsilon_h(t) \Vert_{\infty}\le C_1 \Delta t^{\min\{p,1^*\}}+C_2 h^2;\quad (\Delta t\rightarrow 0, h\rightarrow 0), 
 \end{array} \qquad \Box 
$$
in case of time-dependent BCs.
\end{theorem}

{\sf Proof.} It follows from the PDE-convergence Theorems 4.1 and 4.7 in \cite{GHH-SINUM2020}, which are restricted to the case  of one-stage AMF-W methods. However, by using the same or similar ideas the results can be extended to any $s$-stage AMF-W method with adequate stability properties and a consistency ODE order $p\ge 2$.   \hfill $\Box$

\bigskip
\begin{remark}{\rm It is important to notice that the presence of $\mathcal{O}(h^2)$ in the global error expressions of the above theorems \ref{th-1} and \ref{th-2} comes from the space discretization based on second order central differences. For stable higher order space discretizations, this term would keep the corresponding order.  

It is also remarkable that the result in theorem  \ref{th-2} \  for the case of time independent Dirichlet boundary conditions does not seem to be optimal in the case of the Maximum norm, since in the Table \ref{table1} \ it is observed a PDE convergence of order three for a  third order method in  classical ODE sense. It is still pending to prove a sharper result for the maximum PDE order attainable in  the case of the   maximum norm. We guess that the maximum PDE convergence order is three. This is based on numerical experiments with several methods of high orders and in the theoretical fact that the PDE order in the maximum norm is bounded by the PDE order in the weighted $\ell_2$-norm, which turns out to be $3.25^*$.
} 
\end{remark}

\section{Mitigating the order reduction through interpolants}
In view of theorems \ref{th-1} and  \ref{th-2}, the idea is to transform the PDE problem into a new problem with homogeneous boundary conditions by using some kind of interpolating function. The main part  of the material in this section follows from some lemmas  in \cite[sect. 5]{GHP-JCAM2021}, but  the results here exposed are explained and illustrated in a more concise way. 

Let us consider the general PDE problem,
\begin{equation}\label{eq3-1}
 \begin{array}{l}
u_t=F(\vec{x},t,u, \partial_{x_1} u, \partial_{x_1x_1} u,\ldots,\partial_{x_d} u, \partial_{x_dx_d} u), \quad (\vec{x},t)\in \Omega\times I^*\\[0.5pc]
\gamma(\vec{x},t,u,\partial_{x_1} u,\ldots,\partial_{x_d} u)=\beta(\vec{x},t), \quad (\vec{x},t)\in \partial \Omega\times I^*\\[0.5pc]
u(\vec{x},0)=u_0(\vec{x}), \quad  \vec{x}\in  \Omega=(0,1)^d.
\end{array}
\end{equation} 
 \begin{enumerate}
 \item   Build $\varphi(\vec{x},t) $ in $\bar \Omega\times I^*$ s.t.  $$\gamma(\vec{x},t,\varphi,\partial_{x_1}\varphi,\ldots,\partial_{x_d}\varphi )=\beta(\vec{x},t), \quad (\vec{x},t)\in \partial \Omega\times I^*.$$
\item Define:  $ w(\vec{x},t)=u(\vec{x},t)-\varphi(\vec{x},t)$ on $\bar\Omega\times I^*$
 \end{enumerate}
Then,  assuming the linearity of the boundary conditions in $\gamma$ regarding $u$ and $\partial_{x_j} u,\;j=1,\ldots,d$, we have that 
\begin{equation}\label{eq3-2} \begin{array}{rcl}
w_t&=&F\Big(\vec{x},t,(w+\varphi), \partial_{x_1}(w+\varphi), \partial_{x_1x_1}(w+\varphi),\ldots , \partial_{x_d}(w+\varphi), \partial_{x_dx_d}(w+\varphi)\Big)\\[0.3pc] &-& \partial_t\varphi(\vec{x},t), \end{array}$$  $$ \begin{array}{ll}
\gamma(\vec{x},t,w,\partial_{x_1} w,\ldots,\partial_{x_d} w)=0, & \quad (\vec{x},t)\in \partial \Omega\times I^*, \\[0.5pc]
w(\vec{x},0)=u_0(\vec{x})-\varphi(\vec{x},0), & \quad  \vec{x}\in  \Omega.
\end{array} \end{equation} 
In fact, for the PDE problem (\ref{eq1-0}) we deduce that
\begin{equation}\label{eq3-3} \begin{array}{l}w_t=\mathcal{L} w + r^*(\vec{x},t,w), \quad (\vec{x},t)\in \Omega \times I^*,\quad  \hbox{\rm with }\\[0.3pc] r^*(\vec{x},t,w)=\mathcal{L} \varphi + r(\vec{x},t,w+\varphi)-  \partial_t\varphi,\\[0.3pc]
 \gamma(\vec{x},t,w,\partial_{x_1} w,\ldots,\partial_{x_d} w)=0,  \quad (\vec{x},t)\in \partial \Omega\times I^*,\\[0.4pc]  
 w(\vec{x},0)=u_0(\vec{x})-\varphi(\vec{x},0),  \quad  \vec{x}\in  \Omega.
\end{array} \end{equation}

%%%%%%%%%%%%%%%%%%%%%%%%%%%%%%%%%%%%%%%% 

\medskip

\subsection{How can the interpolant $\varphi$ be built in the interior points of  $\Omega$?}
We consider the case  of Dirichlet, Neumann or Robin conditions  with constant coefficients and the following linear operator to define the boundary conditions, 
$$\begin{array}{rcc} \mathcal{B}&:& C^1(\bar \Omega\times I^*) \rightarrow C(\bar \Omega\times I^*),\\[0.5pc]& & \varphi \rightarrow  \mathcal{B} \varphi \end{array} $$ which is defined on the boundary points $\vec{x}=(x_1,x_2,\ldots,x_d) \in \partial \Omega$ ($x_j=k$ for some $j$, with $k=0$ or $k=1$) as    $$  \begin{array}{l}
\mathcal{B}_k^{(j)}\varphi :=p_k^{(j)}\cdot \varphi\rfloor_{(x_j=k)}+ q_k^{(j)}\cdot \partial_{x_j}\varphi\rfloor_{(x_j=k)}, \; \hbox{\rm where } |p_k^{(j)}|+|q_k^{(j)}|>0,   \\[0.5pc] k=0 \;\mbox{\rm or } k=1 \;\hbox{\rm indicates interval end-points} \; \mbox{\rm and } j=1,2,\ldots,d, \;\hbox{\rm indicates direction.}  
\end{array} $$
The boundary conditions imposed are then given by 
\begin{equation}\label{eq3-3a}\mathcal{B}\varphi(\vec{x},t)=\beta(\vec{x},t), \quad (\vec{x},t)\in \partial \Omega \times I^*. \end{equation}
For instance, the case of Dirichlet or Neumann BCs are respectively given by 
\begin{enumerate}
\item Dirichlet BCs: $p_k^{(j)}=1;\;q_k^{(j)}=0, \;(k=0,1;\;j=1,2,\ldots,d).$
\item Neumann BCs: $p_k^{(j)}=0;\;q_k^{(j)}=1, \;(k=0,1;\;j=1,2,\ldots,d).$
\end{enumerate}

\medskip

\subsection{Building the interpolant}

\smallskip

Based on the following modified result from \cite{GHP-JCAM2021}, we can state 
\begin{theorem}\label{th-3} 
Define, 
\begin{enumerate}
\item $u^{[0]}(\vec{x},t)=\beta(\vec{x},t),\; (\vec{x},t)\in \partial \Omega\times I^*,\quad \beta \in  C^r( \partial \Omega \times I^*).$
\item For $j=1,2\ldots,d$ and $(\vec{x},t)\in \bar \Omega\times I^*$, define:
\begin{enumerate}
\item $\displaystyle \varphi_j(\vec{x},t)=P_j(x_j)	\cdot \mathcal{B}_0^{(j)} u^{[j-1]}(\vec{x},t) + Q_j(x_j) \cdot\mathcal{B}_1^{(j)} u^{[j-1]}(\vec{x},t)$\\[0.5pc]
$\displaystyle u^{[j]}(\vec{x},t)=u^{[j-1]}(\vec{x},t)-\varphi_j(\vec{x},t)$, and 
\item $\displaystyle \varphi(\vec{x},t)=\sum_{j=1}^d \varphi_j(\vec{x},t),$
\end{enumerate}
\item where the one variable polynomials $P_j(\xi), \:Q_j(\xi)$ must fulfil 
$$ \begin{array}{c} p_0^{(j)}P_j(0)+  q_0^{(j)}P'_j(0)=1, \quad  p_0^{(j)}Q_j(0)+  q_0^{(j)}Q'_j(0)=0\\
 p_1^{(j)}P_j(1)+  q_1^{(j)}P'_j(1)=0, \quad  p_1^{(j)}Q_j(1)+  q_1^{(j)}Q'_j(1)=1. \end{array}$$
\end{enumerate}
Then, \ $\varphi(\vec{x},t)\in C^\infty(\Omega\times I^*)\cup C^r(\bar\Omega\times I^*) $ and satisfies \\ \centerline{ $\mathcal{B}\varphi(\vec{x},t)=\beta(\vec{x},t), \quad (\vec{x},t)\in \partial \Omega\times I^*.$}
\end{theorem}

\noindent{\sf Proof.} It follows from Lemmas 2, 3 and 4  in \cite{GHP-JCAM2021}. \hfill $\Box$

\medskip

It is worth to mention that the polynomials in  item 3 of  theorem \ref{th-3} are not unique, but they  always exist and have  a low degree, not greater than three in the worst situation. 

\bigskip
%%%%%%%%%%%%%%%%%%%%%%%%%%%%%%%%%%%%%%%% 

\noindent {\sf A 2D example with mixed BCs of Dirichlet and Neumann type}

\smallskip

The prescribed boundary conditions in the unit square are
$$\begin{array}{ll}
u(0,y,t)=D_0(y,t), \; u(1,y,t)=D_1(y,t),\quad 0\le y\le 1,  & \mbox{\rm (Dirichlet BCs),} \\ 
u_y(x,0,t)=N_0(x,t), \; u_y(x,1,t)=N_1(x,t),\; 0\le x\le 1,  & \mbox{\rm (Neumann BCs).} \end{array}
$$
The interpolant, according the previous result,  is given by 
$$ \varphi(x,y,t)= \sum_{j=1}^2 \varphi_j(x,y,t)$$
with 
$$\begin{array}{rcl}  \varphi_1(x,y,t)&=&(1-x)D_0(y,t) + xD_1(y,t)\\[0.5pc]
\varphi_2(x,y,t)&=&\displaystyle  -\frac{(y-1)^2}{2}\left( N_0(x,t) +(x-1)N_0(0,t)-xN_0(1,t)\right) \\[0.5pc] &+&  \displaystyle \frac{y^2}{2}\left( N_1(x,t) +(x-1)N_1(0,t)-xN_1(1,t)\right) 
 \end{array}
$$

\bigskip

\section{Boundary corrections based on extending the operator to the boundary}\label{sect-4}
Here, we only consider the case of Dirichlet boundary conditions which was shortly considered in  \cite{GHP-NUMAL2023}. The case of  Neumann and Robin boundary conditions are under study  at present.  From (\ref{eq1-1}) it follows that, 
\begin{equation} \label{eq4-1} \begin{array}{l}
\dot V(\vec{x},t)=\mathcal{L}^{(h)} V(\vec{x},t) + r^{(h)}(\vec{x},t,V), \quad (\vec{x},t)\in \Omega_h\times I^*\\[0.5pc]
V(\vec{x},t)=\beta(\vec{x},t), \quad (\vec{x},t)\in \partial \Omega_h\times I^*,\qquad \mbox{\rm (Dirichlet BCs)}\\[0.5pc]
V(\vec{x},0)=u_0(\vec{x}), \quad  \vec{x}\in  \Omega_h. 
\end{array} 
\end{equation} 
It should be observed that this system is not a differential system itself, but it can be cast as an index one Differential Algebraic Equation (DAE), when considering all $V$-entries as unknowns, including those ones on the boundary. By taking into account that many numerical methods (in particular Rosenbrock methods) may suffer from  reduction in the convergence order on DAEs, see e.g. \cite{HW-2}, the situation could  worsen when the DAEs  have variable dimensions and we consider PDE convergence orders, i.e., the situation where $N\simeq (\Delta x_1\Delta x_2 \ldots \Delta x_d)^{-1} \rightarrow \infty$ and the methods are of splitting type. This order reduction was already observed in Table \ref{table2}.

To avoid such an undesirable effect, we propose the following modification which transforms the  DAE problem into an ODE system. This modification is based on empirical results obtained by using several AMF-W-methods on a few 2D,3D, 4D PDEs used as test problems and in some theoretical aspects considered in \cite{GH-APNUM2022, GHH-ESAIM2023} which will be briefly commented below. This modification mainly affects  the treatment of the  boundary conditions,   
\begin{equation} \label{eq4-2}
 \begin{array}{l} 
\dot V(\vec{x},t)=\mathcal{L}^{(h)} V(\vec{x},t) + r^{(h)}(\vec{x},t,V), \quad (\vec{x},t)\in \Omega_h\times I^*\\[0.5pc]
\dot V(\vec{x},t)=\dot \beta(\vec{x},t)+ \Big(\tilde{\mathcal{L}}^{(h)} V(\vec{x},t)-  \tilde{\mathcal{L}}^{(h)} \beta(\vec{x},t)\Big), \quad (\vec{x},t)\in \partial \Omega_h\times I^*,\\[0.5pc]
V(\vec{x},0)=u_0(\vec{x}), \quad  \vec{x}\in  \Omega_h \cup \partial \Omega_h, 
\end{array} 
\end{equation}
where the extension of  $\mathcal{L}^{(h)}$ to the boundary is made as   
\begin{equation}\label{eq4-3} \tilde{\mathcal{L}}^{(h)} V\rfloor_{x_i=0,1 \atop 0<x_j<1,\:j\ne i }= \sum_{j=1\atop j\ne i}^d \Big(\textcolor{black}{a_j(\vec{x},t)}\cdot  \partial^{(h)}_{x_jx_j} V + \textcolor{black}{b_j(\vec{x},t}) \cdot \partial^{(h)}_{x_j} V\Big).
\end{equation}
For other points of the boundary the sum in (\ref{eq4-3}) runs for all indexes except the ones corresponding to components on the boundary, e.g., if the boundary point $\vec{x}=(x_1,x_2,\ldots,x_d)$ has $x_1=0$ and $x_2=0$,  then the above sum in (\ref{eq4-3}) runs from $j=3$ to $d$. 

Observe that we apply differentiation of the boundary conditions and extend the operator to the boundary  in all the allowed directions in a natural way. 

%%%%%%%%%%%%%%%%%%%%%%%%%%%%%%%%%%%%%%%%%%
In this case, after the space discretization  and considering directional splitting, we get the ODE system, 
\begin{equation}\label{eq4-4}
 \displaystyle \dot V(t)=\sum_{j=0}^d F_j(t,V),\quad V(0)=V_0, \quad V(t)\equiv V(\vec{x},t), \;\vec{x} \in \bar \Omega_h,
\end{equation}
where for $j=1,2,\ldots,d$, we perform a directional splitting 
\begin{equation}\label{eq4-5}
\begin{array}{l}
 F_j(t,V)= \left\{ \begin{array}{lr}
\textcolor{black}{a_j(\vec{x},t)} \: \partial^{(h)}_{x_jx_j} V + \textcolor{black}{b_j(\vec{x},t)}  \: \partial^{(h)}_{x_j} V, &  \vec{x}\in\Omega_h\\[0.5pc] \textcolor{black}{a_j(\vec{x},t)} \: \tilde\partial^{(h)}_{x_jx_j} V + \textcolor{black}{b_j(\vec{x},t)}  \: \tilde\partial^{(h)}_{x_j} V ,&\vec{x}\in \partial \Omega_h
\end{array} \right.
\end{array} 
\end{equation}
with 
\begin{equation}\label{eq4-6} \begin{array}{c} \tilde\partial^{(h)}_{x_jx_j} V(\vec{x},t) =\left\{ \begin{array}{ll} 0, & \hbox{\rm if }  x_j\in \{0,1\} \\[0.5pc]  \partial^{(h)}_{x_jx_j} V(\vec{x},t) , & \hbox{\rm otherwise,}  \end{array} \right. \\[1.5pc]
\tilde\partial^{(h)}_{x_j} V(\vec{x},t) =\left\{ \begin{array}{ll} 0, & \hbox{\rm if }  x_j\in \{0,1\} \\[0.5pc]  \partial^{(h)}_{x_j} V(\vec{x},t) , & \hbox{\rm otherwise,}  \end{array} \right. 
\end{array}
\end{equation}
whereas for $F_0$ we consider the reaction term in the interior of the spatial domain and define the following extension of the operator  (denoted $\tilde r^{(h)}$) by using  the Dirichlet boundary conditions
\begin{equation}\label{eq4-7}
\begin{array}{l}
 F_0(\vec{x},t,V)= \left\{ \begin{array}{lr}
r^{(h)}(\vec{x},t,V) &  \vec{x}\in\Omega_h\\[0.5pc] \tilde r^{(h)}(\vec{x},t) &\vec{x}\in \partial \Omega_h,
\end{array} \right.\end{array}
\end{equation}
where  
\begin{equation}\label{eq4-8}\begin{array}{c} \tilde r^{(h)}(\vec{x},t)= \dot \beta(\vec{x},t)-
\tilde{\mathcal{L}}^{(h)} \beta(\vec{x},t), \quad \vec{x}\in \partial \Omega_h. \end{array}
\end{equation}

Now, the time integration is performed with some AMF-W method $(A,L,b,\theta)$  as it is indicated in (\ref{eq1-7}), (\ref{eq1-8}) and (\ref{eq1-8a}), where we use the splitting in (\ref{eq4-4}),  (\ref{eq4-5}), (\ref{eq4-6}), (\ref{eq4-7}) and (\ref{eq4-8}).

We  remark that after an integration step is completed, the AMF-W method is  applied with projection on the boundary to recover the exact values on it, i.e., $V_{n+1}$ is redefined at the boundary points as 
\begin{equation}\label{eq4-11} V_{n+1}=\beta(\vec{x},t_{n+1}), \quad \vec{x} \in \partial \Omega_h.\end{equation}

%%%%%%%%%%%%%%%%%%%%%%%%%%%%%%%%%%%%%%%%%%

\bigskip 

\section{Additional numerical experiments}\label{sect-5}
 
This section is devoted to illustrate the boundary correction technique of the previous Section \ref{sect-4} by showing that it allows to recover  the PDE convergence orders in weighted $\ell_2$-norm and in the maximum norm of AMF-W-methods at the level of time independent boundary conditions. To this aim, we consider  two AMF-W methods widely used in the literature, see e.g. \cite[Chapt. IV]{HV-2003}, \cite{GHP-JCAM2017}, having  classical ODE orders 3 and 4, respectively. We use as test problems three PDE problems which are linear in the diffusion term but with two of them having a nonlinear reaction term.  The PDE problems are semi-discretized in space with a fourth order stencil combined with a second order stencil for the adjacent points to the boundary as it was indicated in section \ref{sect1-4}.
 
\subsection{The PDE problems}
\begin{enumerate}
\item {\sf Problem 1} is the 3D PDE problem given in (\ref{eq1-10}), (\ref{eq1-11}) and (\ref{eq1-12}), by taking $C=1$ (case of time dependent BCs). 
\item {\sf Problem 2}  is a reaction diffusion problem, similar to \cite[p.367]{HV-2003}, with exact solution given by 
$$ u=\beta(x,y,t)=1/(1+\exp(x+y-t)), \quad \Omega=(0,1)^2, \quad t\in I^*=[0,1].$$ 
\begin{equation}\label{eq5-1} \begin{array}{rll}
 u_t(x,y,t)&=&u_{xx}+u_{yy} + u(1-u)(4u-1), \quad (x,y,t)\in \Omega\times I^*\\[0.5pc]
 u(x,y,t)&=&\beta(x,y,t), \quad  (x,y,t)\in \partial \Omega \times I^*\\[0.5pc] 
 u(x,y,0)&=&\beta(x,y,0), \quad (x,y)\in \Omega. 
\end{array}
\end{equation}
\item {\sf Problem 3}  is the 3D version of Problem 2, with exact solution given by 
$$ u=\beta(x,y,z,t)=1/(1+\exp(x+y+z-t)), \quad \Omega=(0,1)^3, \quad t\in I^*=[0,1].$$ 
\begin{equation}\label{eq5-2} \begin{array}{rll}
 u_t(x,y,z,t)&=&u_{xx}+u_{yy} +u_{zz} + 2u(1-u)(3u-1), \quad (x,y,z,t)\in \Omega\times I^*\\[0.5pc]
 u(x,y,z,t)&=&\beta(x,y,z,t), \quad  (x,y,z,t)\in \partial \Omega \times I^*\\[0.5pc] 
 u(x,y,z,0)&=&\beta(x,y,z,0), \quad (x,y,z)\in \Omega. 
\end{array}
\end{equation}
\end{enumerate}

\begin{table}[h!]
\begin{center}
\begin{tabular}{l|l|l} $\Delta t=\Delta x_j$ & \ GE$_{2,h} \quad$ ($p_2$) & \ GE$_\infty \qquad$ ($p_\infty$) \\ \hline 
$1/4 $ & 0.31e+0 \;(---) &  0.51e+0 $\quad $(---) \\
$1/8 $ & 0.58e-1 \;(2.42) &  0.11e+0 $\;$ (2.17) \\
$1/16 $ & 0.95e-2 \;(2.74) &  0.20e$\:$-1 $\;\:$ (2.74) \\$1/32 $ & 0.14e-2 \;(2.74) &  0.29e$\:$-2 $\;\:$ (2.74) \\$1/64 $ & 0.20e-3 \;(2.84) &  0.39e$\:$-3 $\;\:$ (2.91) \\$1/128 $ & 0.27e-4 \;(2.90) &  0.48e$\:$-4 $\;\:$ (3.01) \\ \hline
\end{tabular}
\caption{Case $C=1$ for the 3D linear test problem (\ref{eq1-10})-(\ref{eq1-12}). Global errors  at the end-point and observed convergence orders (in parenthesis) for the corresponding norms with the {\sf AMFW-HV} method.} 
\label{table2-1}
\end{center}
\end{table}

\subsection{The AMF-W methods used}
We are going to consider two AMF-W methods coming from the Rosenbrock family (or W-method family) that have classical orders three and four for ODEs. There are many splitting methods of classical order two, see e.g. \cite[Chapt. IV]{HV-2003}, and most of them preserve the same PDE order of convergence for time-independent boundary conditions, but splitting methods with time PDE order higher than two are scarse in the literature.

\begin{enumerate}

\item {\sf AMFW-HV} is the two-stage AMF-W $(A,L,b,\theta)$ method of classical order 3, considered in (\ref{eq1-9}), see also \cite[p. 400]{HV-2003}. It has PDE order three in the weighted $\ell_2$-norm \cite{GHH-ESAIM2023}, see also the Theorem \ref{th-1} in the current paper. In the maximum norm it seems to have also PDE order three in time, but only order two can be guaranteed from Theorem \ref{th-2}. It is a pending question to show order three as the numerical experiments seem to indicate.  
\item {\sf AMFW-3/8} is the four stage AMF-W $(A,L,b,\theta)$ method of classical order 4 based in the 3/8 Runge-Kutta given in  \cite[p. 153-154, (35), (37) and (38)]{GHP-JCAM2017}. The method has order four as Rosenbrock method and order three as W-method \cite{STEIHAUG-WOLFB-1979} when arbitrary $W$ matrices are used to replace the current Jacobian $\partial_V F(t_n,V_n)$. The selected value of $\theta$ and the remaining coefficient matrices are
 \begin{equation} \label{eq5-5} \begin{array}{c}   \displaystyle{ \theta =\frac{1}{2}}\;, \qquad  b^\top= \left( 7/8, \:9/8,\: 9/8,\: 1/8\right),  \\[0.7pc]
 A=\left(\begin{array}{cccc} 0 & 0 & 0 &\\ 4/9 & 0& 0 &0 \\ -1/3&1 & 0 & 0\\ -1 &-3 & 6& 0 \end{array} \right), \quad  L=\left(\begin{array}{cccc} 0 & 0 & 0 &0 \\ -4/3 & 0 & 0 &0 \\ -1 &-1 & 0 & 0\\ 1 & -3 & -6 & 0 \end{array} \right). 
 \end{array} \end{equation}
\end{enumerate}

This method has PDE order $3.25^*$ in the weighted $\ell_2$-norm as it is deduced from Theorem \ref{th-1} and from the  analysis carried out in \cite[sect. 5-6]{GHP-JCAM2017}. According to Theorem \ref{th-1}, it has PDE order two   at least in the maximum norm. However, it appears to have order three for 2D and 3D  problems.

\begin{table}[h!]
\begin{center}
\begin{tabular}{l|l|l} $\Delta t=\Delta x_j$ & \ GE$_{2,h} \quad$ ($p_2$) & \ GE$_\infty \qquad$ ($p_\infty$) \\ \hline 
$1/4 $ & 0.44e-1 $\;$ (---) &  0.66e-1 $\quad \; $(---) \\
$1/8 $ & 0.44e-2 \;(3.32) &  0.77e-2 $\;\;\:$ (3.01) \\
$1/16 $ & 0.58e-3 \;(2.93) &  0.12e$\:$-2 $\;\:$ (2.70) \\$1/32 $ & 0.79e-4 \;(2.88) &  0.15e$\:$-3 $\;\:$ (2.96) \\$1/64 $ & 0.84e-5 \;(3.23) &  0.12e$\:$-4 $\;\:$ (3.69) \\$1/128 $ & 0.79e-6 \;(3.41) &  0.11e$\:$-5 $\;\:$ (3.48) \\ \hline
\end{tabular}
\caption{Case $C=1$ for the 3D linear test problem (\ref{eq1-10})-(\ref{eq1-12}). Global errors  at the end-point and observed convergence orders (in parenthesis) for the corresponding norms with the {\sf AMFW-3/8} method.} 
\label{table3-1}
\end{center}
\end{table}

\medskip

\medskip

\subsection{Estimating  the spatial errors of the nonlinear PDE problems}
We expect that the spatial errors   are of order  four due to the local  order four of the spatial discretization made. To assess this fact, we will assume that the spatial errors (SE) of the ODEs regarding the PDE behave as
\begin{equation}\label{eq5-10}
SE(t^*,h):=u(\vec{x}_G,t^*)-V(\vec{x}_G,t^*)= \Theta(\vec{x}_G,t^*)h^p+ \mathcal{O}(h)^{p+1}, \quad p\ge 2.
\end{equation}
with $\Theta(\vec{x}_G,t^*)$ being a vector with components coming from the grid-values taken for  some smooth function $\Theta(\vec{x},t)$ defined on $\bar \Omega\times I^*$. Then, to compute the order $p$ numerically, we integrate the space discretized ODEs with  the  AMFW-3/8 method. We assume for the global errors (spatial errors plus temporal errors)  that they behave   (in the $\ell_2$-norm) as 
\begin{equation}\label{eq5-11}
GE(t^*,h,\Delta t):=u(\vec{x}_G,t^*)-V(\vec{x}_G,t^*)= \Theta(\vec{x}_G,t^*)h^p+ \mathcal{O}(h)^{p+1} + \mathcal{O}(\Delta t)^{3}, \quad h,\:\Delta t \rightarrow 0,
\end{equation}
where the vector $\Theta(\vec{x}_G,t)$ is as indicated in (\ref{eq5-10}) and the term $\mathcal{O}((\Delta t)^3)$, comes from the time integration with the {\sf AMFW-3/8} rule as we have seen in the time integration of the  3D linear problem. By performing  time integrations with a time step-size $\Delta t=\kappa\cdot h^{5/3}$ ($\kappa=1,0.5, 0.25$, or some other constant), we get that $(\Delta t)^3=\mathcal{O}(h^{5})$, and from (\ref{eq5-11}) it follows that   
\begin{equation}\label{eq5-12}
GE(t^*,h,\kappa h^{5/3}):= \Theta(\vec{x}_G,t^*)h^p+ \mathcal{O}(h^{q}), \;q=\min\{p+1,5\},\quad h\rightarrow 0.
\end{equation} 
 In such a situation for $p<5$, we compute the  order of convergence of the space discretizations by  
\begin{equation}\label{eq5-13}
p\simeq \log (\|GE(t^*,2h)\|/ \|GE(t^*,h)\|) /\log 2.   
\end{equation}

\begin{table}[h!]
\begin{center}

\begin{tabular}{l|lc|lc} $h=\Delta x_j$ & \ GE$_{2,h}$ & ($p_2$) & \ GE$_\infty$ &($p_\infty$) \\ \hline 
   1/8 &   0.1372D-05 & (---) &   0.3163D-05     & (---) \\      
    1/16 &  0.9562D-07 & (3.84) &  0.2102D-06    & (3.91) \\
    1/32 &  0.6311D-08 &(3.92) &    0.1401D-07    &    (3.91)   \\
    1/64 &  0.4066D-09&  (3.96) & 0.8920D-09    &(3.97)\\
    1/128 &  0.2585D-10 & (3.98) &   0.5616D-10    & (3.99)
\\    
\hline
\end{tabular}

\caption{For Problem 2, we estimate at the end-point $t^*=1$, the global errors (GE) for the space discretization and  the estimated convergence orders $p$  for  $h=\Delta x_j,\:j=1,2$ and $\Delta t=h^{5/3}$, computed according to (\ref{eq5-12})-(\ref{eq5-13}).}  
\label{tablespace-1}
\end{center}
\end{table}

\begin{table}[h!]
\begin{center}

\begin{tabular}{l|lc|lc} $h=\Delta x_j$ & \ GE$_{2,h}$ & ($p_2$) & \ GE$_\infty$ &($p_\infty$) \\ \hline

   1/8 &   0.1786D-05& (---)&    0.3382D-05   & (---)\\
  1/16 &   0.1167D-06 & (3.94)&    0.2215D-06&       (3.93)\\
  1/32 &  0.7476D-08 &(3.96)&   0.1449D-07  &  (3.94)\\
  1/64 & 0.4744D-09 & (3.98)&    0.9286D-09   & (3.96) \\ 
  1/128 & 0.2991D-10 & (3.99) &  0.5940D-10   &    (3.97)\\    
\hline
\end{tabular}

\caption{For Problem 3, we estimate at the end-point $t^*=1$, the global errors (GE)  for the space discretization and the estimated convergence orders $p$  for  $h=\Delta x_j,\:j=1,2,3$ and $\Delta t=\frac{1}{4} h^{5/3}$, computed according to (\ref{eq5-12})-(\ref{eq5-13}).}  
\label{tablespace-2}
\end{center}
\end{table}

From Table \ref{tablespace-1} and Table \ref{tablespace-2},  it is appreciated that the  space discretization errors are of order four in both norms, the weighted $\ell_2$-norm and  the maximum norm, independently of the number of space variables (2D or 3D). This kind of space discretization is quite nice and advantageous over the standard central second order one since, although the computational effort involved with it is slightly more expensive per  linear system solution than the one based on standard second order differences, only systems with banded pentadiagonal matrices  instead of tridiagonal banded matrices need to be solved for the fourth order discretization, while having the  advantage that the new one seems to gain two orders in the spatial order of convergence. It is worth to mention that a rigorous proof of the fourth order of convergence for the new space discretization is still in progress.

\subsection{Estimating the temporal order of the global errors in PDE sense}

Due to the fact that the AMF-W methods applied with boundary correction based on interpolants yield  similar global errors as in the case of homogenous boundary conditions and their orders  for linear problems are theoretically well stablished in \cite{GHH-ESAIM2023} (see also here, Theorem \ref{th-1} and Theorem \ref{th-2}), we will  not pursue this issue in our experiments. We will illustrate the case of the AMF-W  methods modified with the boundary corrections proposed in section \ref{sect-4}. 
   
According to the results from the previous section and from  section \ref{sect-2}  (case of the linear problem) we will assume that the global errors of the MoL approach (space discretization plus time discretization) behave as, see (\ref{eq5-11}), 
\begin{equation}\label{eq6-1} \begin{array}{l}
GE(t^*,h,\Delta t):=u(\vec{x}_G,t^*)-V(\vec{x}_G,t^*)= \Theta(\vec{x}_G,t^*)h^4(1+ \mathcal{O}(h))+ \\[0.5pc] \qquad \qquad \qquad \qquad\chi_V(\vec{x}_G,t^*) \Delta t^p(1+\mathcal{O}(\Delta t^{r})), \;r>0,\end{array}
\end{equation}
with $\Theta(\vec{x}_G,t^*)$ and $\chi_V(\vec{x}_G,t^*)$ being vectors with components coming respectively from smooth functions defined on $\bar \Omega\times I^*$. We can then estimate the temporal order of convergence $p$ from formula (\ref{eq1-16}) by using $$\Delta t= \kappa h, \; h=\Delta x_j, \;j=1,2,\ldots,d, \quad \kappa =2,1, 0.5 ,0.25, \ldots$$

For the linear problem 1 the term $\Theta(\vec{x}_G,t^*)$ vanishes due to the fact that there are no spatial errors in the discretizations. From Table \ref{table2-1}  we can appreciate that the {\sf AMFW-HV} method behaves with order three in both norms as it was expected. From Table \ref{table3-1}  we can appreciate that the {\sf AMFW-3/8} method behaves with some erratic orders among 2.9 and 3.5, but the trend is towards orders greater than 3. It is expected order $3.25^*$ in the weighted $\ell_2$-norm and order 3 in the maximum norm. This can be better appreciated for problem 2 and method  {\sf AMFW-3/8} in Table \ref{table3-2}. For problem 2 it can be also seen  in Table \ref{table2-2} that  the {\sf AMFW-HV} method  keeps orders tending towards 3 in both norms as it was expected. 

However, in problem 3 (3D nonlinear case) the temporal orders cannot be appreciated in  Table \ref{table2-3} and Table \ref{table3-3} for any of both AMF-W methods as it was in the case of problem 2. Due to storage limitations and the time taken for the long computations in problem 3 when $h$ is decreased beyond $h_0=1/230$ (this case involves  ODE systems of dimension $N_0\simeq 12\cdot 10^6$ equations) and due to   the uncertainty about the temporal orders  on this problem displayed in Tables \ref{table2-3} and \ref{table3-3}, which are  based on formula (\ref{eq1-16}), we have designed another way of estimating the temporal orders of convergence based on the global error formula in   (\ref{eq6-1}). To estimate the temporal order $p$, we consider for each $h=\Delta x_j, \;j=1,2,3$, three time integrations with time stepsizes $\Delta t_j=2^{j-1}h, \;j=0,1,2$. Then from  (\ref{eq6-1}) we can deduce that 
\begin{equation}\label{eq6-2} \begin{array}{c}
p\simeq \log (\|\epsilon(t^*,h,h)\|/ \|\epsilon(t^*,h,h/2)\|) /\log 2,\\[0.5pc]
\epsilon(t^*,h,\Delta t)=V(\vec{x}_G,t^*,h,2\Delta t)-
V(\vec{x}_G,t^*,h,\Delta t). \end{array}
\end{equation}

\begin{table}[h!]
\begin{center}
\begin{tabular}{l|lc|lc} $\Delta t=\Delta x_j$ & \ GE$_{2,h}$ & ($p_2$) & \ GE$_\infty$ &($p_\infty$) \\ \hline

     1/8  & 0.1001D-04& (---)&  0.2695D-04 & (---)\\  
    1/16 &  0.1184D-05 &(3.08)&   0.3326D-05&       (3.02)\\ 
    1/32  & 0.1383D-06 &(3.10)&   0.4067D-06   &    (3.03)\\  
       1/64 &  0.1715D-07 & (3.01)&  0.5105D-07&  (2.99)\\
   1/128  & 0.2258D-08& (2.93)&  0.7404D-08& (2.79) \\
   1/256  & 0.3048D-09 &(2.89) &  0.1074D-08  &     (2.79)\\
   1/512 &  0.4108D-10 &(2.89) &   0.1528D-09    &(2.81) \\ 
  1/1024  & 0.5472D-11 & (2.91)&  0.2125D-10 & (2.85)\\
\hline
\end{tabular}
\caption{Global errors  at the end-point  and observed convergence orders (in parenthesis) for Problem 2, with the {\sf AMFW-HV} method for $\Delta t=h=\Delta x_j,\:j=1,2$. } 
\label{table2-2}
\end{center}
\end{table}

\begin{table}[h!]
\begin{center}
\begin{tabular}{l|lc|lc} $\Delta t=\Delta x_j$ & \ GE$_{2,h}$ & ($p_2$) & \ GE$_\infty$ &($p_\infty$) \\ \hline 
    $ 1/8$   &0.29e-05& (---)&   0.54e-05& (---)\\ 
   $ 1/16$ &  0.27e-06&  (3.41)&  0.55e-06& (3.29)\\    
   $ 1/32$&   0.22e-07& (3.61)&   0.57e-07& (3.29)\\    
   $ 1/64$&   0.18e-08&  (3.61)&  0.63e-08& (3.17)\\    
  $ 1/128$&   0.16e-09 &  (3.47)& 0.75e-09& (3.07)\\   
  $ 1/256$& 0.16e-10& (3.34)&  0.93e-10&  (3.01) \\
  $ 1/512$&   0.17e-11&  (3.28)& 0.18e-10& (2.99)\\   
 $ 1/1024$  & 0.17e-12&  (3.26)& 0.15e-11& (2.98)\\  
\hline
\end{tabular}
\caption{Global errors  at the end-point  and observed convergence orders (in parenthesis) for Problem 2, with the {\sf AMFW-3/8 } method with $\Delta t=h=\Delta x_j,\:j=1,2$.} 
\label{table3-2}
\end{center}
\end{table}

\begin{table}[h!]
\begin{center}
\begin{tabular}{l|lc|lc} $\Delta t=\Delta x_j$ & \ GE$_{2,h}$ & ($p_2$) & \ GE$_\infty$ &($p_\infty$) \\ \hline 

    $ 1/8 $&  0.29e-04& (---)&  0.88e-04& (---)\\          
  $  1/16 $&  0.30e-05&  (3.24)& 0.14e-04& (2.60)\\   
    $1/32$&   0.34e-06&  (3.17)& 0.24e-05&(2.61)\\   
    $1/64$&   0.45e-07 &  (2.89)&0.43e-06&(2.48)\\  
  $ 1/128$&   0.65e-08& (2.80)&  0.78e-07&(2.45)\\   
\hline
\end{tabular}
\caption{Global errors  at the end-point  and observed convergence orders (in parenthesis) estimated from (\ref{eq1-16}) for Problem 3, with the {\sf  AMFW-HV} method  with $\Delta t=h=\Delta x_j,\:j=1,2,3$.} 
\label{table2-3}
\end{center}
\end{table}
%-----------------------------
\begin{table}[h!]
\begin{center}
\begin{tabular}{l|lc|lc} $\Delta t=\Delta x_j$ & \ GE$_{2,h}$ & ($p_2$) & \ GE$_\infty$ &($p_\infty$) \\ \hline 
    $ 1/8$&   0.60e-05& (---)&   0.12e-04& (---)\\    
   $ 1/16$&   0.49e-06& (3.60)&  0.16e-05& (2.88)    \\
    $1/32$&   0.30e-07 &(4.05)&  0.25e-06& (2.70)\\
    $1/64$&   0.22e-08& (3.74)&  0.43e-07& (2.54)\\   
  $ 1/128$&   0.24e-09& (3.22)&  0.77e-08&     (2.48) \\  
\hline
\end{tabular}
\caption{Global errors  at the end-point  and observed convergence orders (in parenthesis) estimated from (\ref{eq1-16}) for Problem 3, with the {\sf AMFW-3/8} method and $\Delta t=h=\Delta x_j,\:j=1,2,3$.} 
\label{table3-3}
\end{center}
\end{table}

% Estimating the temporla global errors and its order in pde sense in other way 

\begin{table}[h!]
\begin{center}
\begin{tabular}{l|lc|lc} $h=\Delta x_j$ & \ GE$_{2,h}$ & ($p_2$) & \ GE$_\infty$ &($p_\infty$) \\ \hline 
     1/8 &   0.2862D-04& (3.00) &  0.8784D-04    &  (2.29)\\
    1/16 &  0.3026D-05 & (3.53) &  0.1446D-04    &   (2.56)\\ 
    1/32  & 0.3353D-06 & (3.33) &  0.2372D-05    &(2.70)\\
    1/64 &  0.4515D-07& (3.01) &   0.4267D-06    &  (2.75)\\    
   1/128 &  0.6490D-08 & (2.85)&  0.7802D-07 &      (2.72) \\   
   1/200 &  0.1847D-08 & (2.83)&   0.2813D-07    &  (2.78) \\
      1/224  & 0.1339D-08& (2.84) &   0.2198D-07 &   (2.79)\\
\hline
\end{tabular}
\caption{Global errors  at the end-point  and observed temporal orders of convergence (in parenthesis) for Problem 3, with the {\sf AMFW-HF} method for  $h=\Delta x_j,\:j=1,2,3$ and using $\Delta t=2h,h,h/2$, according to formula (\ref{eq6-2}).} 
\label{table2-4}
\end{center}
\end{table}

\begin{table}[h!]
\begin{center}
\begin{tabular}{l|lc|lc} $h=\Delta x_j$ & \ GE$_{2,h}$ & ($p_2$) & \ GE$_\infty$ &($p_\infty$) \\ \hline 
     1/8 &  0.5953D-05 &(2.58)&  0.1188D-04    &  (2.33)  \\

 1/16   &0.4910D-06& (3.24)&   0.1615D-05  & (2.69)\\  
    1/32 &  0.2970D-07 &(3.92)&  0.2488D-06    &  (2.28)\\  
    1/64 &  0.2219D-08 &(3.46) &   0.4291D-07       &(2.35)\\ 
   1/128 &  0.2386D-09 &(3.14) &  0.7707D-08    &  (2.62)\\ 
   1/200  & 0.5960D-10& (3.13)&   0.2561D-08    & (2.86)\\
 1/224 & 0.4183D-10   & (3.14)&  0.1936D-08  &    (2.93)\\ 
\hline
\end{tabular}
\caption{Global errors  at the end-point  and observed temporal orders of convergence (in parenthesis) for Problem 3, with the {\sf AMFW-3/8} method for  $h=\Delta x_j,\:j=1,2,3$ and using $\Delta t=2h,h,h/2$, according to formula (\ref{eq6-2}).} 
\label{table3-4}
\end{center}
\end{table}

Based on the formula (\ref{eq6-2}) we have displayed in Table \ref{table2-4} and Table  \ref{table3-4} the global errors and the temporal order of convergence estimated for the {\sf AMFW-HV} and {\sf AMFW-3/8} methods, respectively. There, it can be appreciated a trend to the order three for the {\sf AMFW-HV}  method  in both norms when $h$ is decreased. For the case of the  {\sf AMFW-3/8} method, the  order seems to approach to $3.2$ in the weighted $\ell_2$-norm and tends to order $3$ in the maximum norm. It should be observed that now is not needed to halve the space mesh-grid in the integrations to estimate the temporal order due to the fact that it  is computed for the same $h$ over three time  integrations halving the time-step in each one. We have included the  previous entries for $h$ in   Table \ref{table2-4} and Table \ref{table3-4}  (cases $h=1/8, 1/16, 1/32, 1/64, 1/128$) in order to be consistent with the results presented in the previous tables and also to appreciate the changes produced with the new estimation.  It should be observed that for smaller $h$,  more precise  order estimations are collected in Table \ref{table2-4} and Table \ref{table3-4}.

\section{Conclusions, remarks and future research}
We have presented two techniques that completely avoid or mitigate the order reduction on the convergence orders in ODE sense at the level of time independent boundary conditions (BCs) the convergence orders of a Mol approach (Method of Lines) applied on multidimensional parabolic PDEs when time dependent BCs are imposed. According to the observed results, these approaches  present  fourth order of convergence in space and third order of convergence in time at least.  The time integration is based on splitting methods of  AMF-W type with order greater than two. Typically most  of the best methods of splitting type in the literature present order two. The temporal convergence orders are illustrated on two AMF-W-methods and they achieve order three in the maximum norm  and  order $3.25^*$ in the weighted Euclidean norm for the second method (AMFW-3/8 rule). The results are illustrated on several 2D and 3D problems, where two of them present a smooth nonlinear reaction term. We also collect a few theoretical results  supporting some of the convergence orders here observed. One novelty is that the results apparently extend to the new space discretization which seems to be globally of fourth order despite that it is based on a local combination of second order and fourth order central differences. The new space discretization  is  not very demanding computationally and it allows to get  relative high accuracies, namely accuracy $\epsilon\simeq \mathcal{O}(h^4)$,  without using too many space grid-points, which  significantly reduces the computational costs of the time integrations due to the reduction of the dimension of the ODE systems involved in it compared  with those using space discretizations based on central second order differences. 

Research in progress is the extension of the second boundary correction  here proposed (extension of the operator to the boundary) in order to recover the convergence orders when time dependent  Neumann and Robin conditions are imposed. Also some theoretical justification of the fourth order space discretization of the PDE is in progress.   Applications of this MoL approach, extending the fourth order space discretization to some problems with mixed derivatives, frequently appearing in Finance  such as the Heston model and other models with more that two "space" dimensions ($d>2$), in combination with AMF-W-methods are also in mind (related works on this topic are \cite{HOUTW-JCAM2016,LPV-SISC-2021}). 

Another important remark is that both  boundary corrections here proposed apparently work for most of splitting methods as long as they do not suffer order reduction when time independent boundary conditions are imposed, such as the Douglas method, ADI methods,  Strang splitting, Trapezoidal splitting, etc, see e.g.  \cite[Chapt.IV]{HV-2003} and \cite{MARCHUK-1990}. 

%%%%%%%%%%%%%%%%%%%%%%%%%%%%%%%%%%%

\end{document}